\documentclass[10pt,conference]{IEEEtran}
\usepackage{amsmath,amssymb,amsfonts} 
\usepackage{graphicx}                 
\usepackage{verbatim} 
\usepackage{flushend} 
\usepackage{multirow} 
\usepackage{theorem}
\newtheorem{Theorem}{Theorem}

\begin{document}
\title{Demand Shaping to Achieve Steady Electricity Consumption with Load Balancing in a Smart Grid}

\author{\IEEEauthorblockN{Naveed Ul Hassan$^{\dag}$, Xiumin Wang$^{\ddag}$, Shisheng Huang$^{\ddag}$ and Chau Yuen$^{\ddag}$}
\IEEEauthorblockA{$^{\dag}$ Department of Electrical Engineering, SSE, LUMS, Lahore \\
$^{\ddag}$Singapore University of Technology and Design, Singapore
} }
\maketitle

\begin{abstract}
The purpose of this paper is to study conflicting objectives between
the grid operator and consumers in a future smart grid.
Traditionally, customers in electricity grids have different demand
profiles and it is generally assumed that the grid has to match and
satisfy the demand profiles of all its users. However, for system
operators and electricity producers, it is usually most desirable,
convenient and cost effective to keep electricity production at a
constant rate. The temporal variability of electricity demand forces
power generators, especially load following and peaking plants to
constantly manipulate electricity production away from a steady
operating point. These deviations from the steady operating point
usually impose additional costs to the system. In this work, we
assume that the grid may propose certain incentives to customers who
are willing to be flexible with their demand profiles which can aid
in the allowance of generating plant to operate at a steady state.
In this paper we aim to compare the tradeoffs that may occur between
these two stakeholders. From the customers' perspectives, adhering
to the proposed scheduling scheme might lead to some inconvenience.
We thus quantify the customers inconvenience versus the deviations
from an optimal set by the grid. Finally we try to investigate the
trade-off between a grid load balancing objective and the customers'
preferences.

\end{abstract}

\section{Introduction}
Electricity demand in the residential sector can be decomposed into
a combination of individual appliances aggregated by individual
households. These appliances are tied together through different
activities performed by users throughout a day and each of these
activities may involve one or more of these power consuming devices.
These appliances are conventionally managed by each user according
to his/her preferences, e.g. one may decide to wash clothes early in
the morning before he leaves for work, and washing clothes is an
activity or task which involves the use of washing machine, dryer
etc. Different users can perform this task at different hours of the
day according to their convenience. And many of such acvitives/task
are flexible and can be performed at any time during a day. On the
other hand, there may be certain activities which can be regarded as
essential and which needs to be performed daily at exactly specified
time slots e.g. after sunset from 7 pm till mid night one has to
turn on the lights. Such activities and the devices involved in
these activities then contribute towards electricity load which is
essential and which has strict scheduling requirements.

In a traditional grid, the dominant setup has been to serve the
preferences of the users as the priority need and match electricity
supply to the instantaneous demand. This however requires constant
manipulation of electricity production levels. As a consequence
power generating plants suffer large deviations from their steady
operating points which impose additional costs to the overall
system. All this is changing as the grid is becoming smart
\cite{smart_g1}-\cite{smart_gn}. A smart grid can help the operator
in shaping the demand (e.g. schedule the washing machine at a later
time slot when there is less demand) so as to reduce the overall
societal cost for them, this can be done through the flattening of
the demand curve \cite{flatten_1}. To achieve a flatter demand
curve,  it can propose incentives (e.g. discount) to users to change
their preference levels for different activities. Users can then
allow the grid to manage and schedule certain appliances to enjoy
these benefits at the expense of suffering some level of
inconvenience.

In this paper, we attempt to quantify the inconvenience levels,  by
varying the number of appliances that participate through deviation
from their preferred scheduling time slots and also by varying the
number of time slots each activity deviates. We can thus identify a
compromise between the grid operator objectives and user convenience
levels. We believe such understanding is beneficial for the grid to
design effective incentive to achieve load balancing in a smart
grid.

There are some recent studies on this problem. In \cite{study_1}
authors design incentives and propose scheduling algorithms
considering strictly convex functions of costs. Users are given
incentives to move to off peak hours and these incentives are
proposed using game theoretic analysis. However they do not consider
or quantify the inconvenience levels of the users. In \cite{study_2}
authors propose pricing scheme for users in order to achieve a
perfectly flat demand curve. They show that finding an optimum
schedule is NP-hard problem. They propose centralized and
distributed algorithms depending on the degree of knowledge of the
state of the network. The authors in \cite{study_3} propose a
strategy to achieve a uniform power consumption over time. Their
algorithm schedules the devices in such a way that a target power
level is not exceeded in each time slot. However again the authors
do not take into account the inconvenience level of users while
designing these algorithms. In \cite{study_4} the authors use convex
optimization tools and solve a cooperative scheduling problem in a
smart grid. The authors in \cite{study_5} use a water-filling based
scheduling algorithm to obtain a flat demand curve. The proposed
algorithm does not require any communication between scheduling
nodes. The authors also study the possible errors in demand forecast
and incentives for customer participations. It should be noted that
the objective of all these studies is to achieve flat demand curve
for the grid. However in this paper we study the compromise between
the grid objective of flat demand vis-a-vis the user inconvenience
levels, as the acceptance from the users is the key to have smart
grid to be succeed.

The rest of the paper is organized as follows. In section II we describe the
load model, our approach and problem formulation. Proposed solution,
algorithms and metric for comparing various schedules are described in section III.
Simulation results are presented in section IV while the paper is concluded in section V.

\section{Load Model and Problem Formulation}
\subsection{Load Model}
In this paper we consider two types of loads in the grid i.e.
essential and flexible. Essential load is due to essential
activities and the devices involved in these activities have fixed
scheduling needs. Flexible load is due to flexible activities and
the devices involved in these activities can have flexible
scheduling requirements. There is a preferred scheduling time slot
for these flexible activities and user feels most convenient if
these activities are performed according to their preferences.
However we assume a generalized framework that if some activity or
task is declared as flexible then it can be scheduled in time slots
either before or after the preferred time slot for this activity.
For example, pre-cooling a room is an activity that can be scheduled
before the preferred time slot, while cloth washing is an activity
that can be scheduled after the preferred time slot. We understand
that there is no activity that can be scheduled both before and
after the preferred time slot, but in this study, we just assume a
generic load with such flexibility to facilitate the problem
formulation.

The level of inconvenience is measured by the deviation of an
activity from its specified time slot. The more an activity is
scheduled beyond its specified preferred time slot (either to the
left or to the right of it) the more inconvenience a user faces. In
the rest of the paper, the terms, devices, activities and tasks are
used interchangeably. Similarly the terms, flexible and shiftable
are also used interchangeably.

Given a set of tasks and their energy consumptions, we propose two extreme schedules to serve as bounds. The first
schedule is optimal for the grid in terms of load balancing and the second schedule is best for the
user in terms of its preference for non-essential tasks:
\begin{itemize}
    \item \textbf{Grid Convenient (GC) Schedule}: For the given set of essential and shiftable loads, this represents the best schedule from the perspective of the grid. This schedule does not care about the user preferences in scheduling essential as well as non-essential tasks. Instead the objective of this schedule is to achieve maximum load balancing across various time slots. We can obtain this schedule by equally dividing all the load in each time slot.
    \item \textbf{User Convenient (UC) Schedule}: This schedule is the best schedule from the customer's perspective. This is another extreme schedule which does not take into account the load balancing preferences of the grid. Instead it schedules all the non-essential tasks at the most preferred time slots specified by the users. This schedule is most convenient for the users.
\end{itemize}
Any other schedule for the given set of loads will lie between these two extremes. For a given set of essential and shiftable loads, the GC schedule is practically impossible to achieve because there is in reality, not much flexibility in shifting the essential loads. Since we assume that we can only shift the non-essential loads, we study the region between these two extreme schedules through the following parameters:
\begin{itemize}
    \item We change the allowable time slot deviation of non-essential devices from their preferred time slots, serving as a proxy to changing the convenience levels of users. It allows for us to schedule a device within a flexible number of time slots either to the left or to the right of its preferred time slot.
    \item We vary the number of non-essential devices willing to be flexible. All the devices which declare themselves as non-flexible will then be treated as essential loads and will start exactly at their preferred time slots.

\end{itemize}
Through this study, results can influence the stakeholders involved in this system.
The grid can define incentives by measuring the deviation of a given schedule from the perfectly
flat demand profile while also keeping in view the GC schedule for given load conditions.
Similarly a customer can through feedback from its deviation of a given
schedule from the UC schedule, readjust its preference conditions.

\subsection{Problem Formulation}
Let $\mathcal{A}$ denote the set of all essential tasks. We assume that the electricity consumption data
of these essential tasks on an hourly basis are known. Let $E(t) \:,\forall t=1,\ldots,T$ denote the
consumption of electricity by all the essential tasks to be performed during the $t^{th}$ time slot (maybe hour or half hour etc).
Let $\mathcal{S}$ denote the set of all $K$ non-essential tasks.
The electricity consumption of these non-essential tasks is also assumed to be known.
For a non-essential task $i \in \mathcal{S}$, let $\tilde{S}_i$ denote
its total energy consumption. Let $1\leq \Delta T_i \leq T$
denote the total time required to complete non-essential task $i$.
We allow for non-essential tasks to require several time slots to complete,
and once the we decide to carry out this task at time $t$ then we cannot stop it until it is completed.
Let $B_i$ denotes the best operating time for task $i$.
Since we have to finish all the non-essential tasks within $T$ time slots,
therefore we assume that $B_i\leq T-\Delta T_i+1\:,\forall i$ (to allow task $i$ to finish by time $T$).
Let $S_i(t)$ denote the portion of non-essential load $i$ scheduled at time $t$.
Similarly, let $\bar{S}_i=\{S_i(1),\ldots,S_i(T) \}$ contain the per time slot load of
non-essential device $i$. It should be noted that if device $i$ is schedule in time slot
$J$ then,
\begin{eqnarray}
S_i(t)=\Bigg\{\begin{array}{c}
                \frac{\tilde{S}_i}{\Delta T_i},\quad \text{if} \: t=J,\ldots,\Delta T_i-1 \\
                0,\quad \text{otherwise}
                \end{array}
\label{subruly1}
\end{eqnarray}

\subsection{Extreme Schedules}
\subsubsection{Grid Convenient Schedule}
The objective of this schedule is to achieve perfect load balancing for the grid.
This schedule re-distributes the essential as well as flexible load equally in all
time slots. Let us denote the perfectly flat schedule by $\hat{\mathcal{R}}$. It can be obtained as follows:
\[L(t)=\frac{\sum_{t=1}^T E(t)+ \sum_{i\in \mathcal{S}} \tilde{S}_i}{T}\:,\forall t \]
Once again note that this schedule is not a practical schedule for
the given set of essential and shiftable loads. However this
schedule represents the ideal situation for the grid, and merely
serve as benchmark purposes.

\subsubsection{User Convenient Schedule}
The objective of this schedule is to carry out all the essential and
non-essential tasks at their specified best time slots. This
schedule can be determined by treating the non-essential tasks like
essential load at the specified time slots. E.g. if for task $i$ the
best time slot is $B_i=3$ and $\Delta T_i=2$ then,
\[\bar{S}_i^e=\{0,0,\tilde{S}_i/2,\tilde{S}_i/2,0,0,\ldots,0 \} \]
Let us denote this schedule by $\tilde{\mathcal{R}}$.
We determine $\bar{S}_i^e=\{S_i(t)\}_{t=1}^{T} \:,\forall i$ and then the total scheduled load during time slot $t$ is given as,
\[L(t)=E(t)+ \sum_{i \in \mathcal{S}} S_i(t) \quad, \forall t \]
This is a practical schedule, representing the current status quo and the most convenience for the users.

\subsection{Practical Schedules}
We can obtain a range of schedules between the above two extreme schedules by
changing the number of devices declaring themself as flexible and also by defining
the number of time slot deviations they are willing to tolerate. If all the devices
declare them as non-flexible then we will obtain schedule $\tilde{\mathcal{R}}$ (UC Schedule).
On the other hand if all the all non essential devices declare them as flexible and are willing to
tolerate maximum possible time slot deviation then such a schedule, though not
perfectly flat (due to the presence of essential loads in each time slot) will
be the best schedule for the grid for a given set of loads.
Let $\hat{\mathcal{S}} \subseteq \mathcal{S}$ denote the set of devices
which declare themselves as flexible. Similarly let $X_i$ denote the time slot
deviation that device $i \in \hat{\mathcal{S}}$ is willing to tolerate.
It means that we aim to schedule non-essential task $i$ within $X_i$ time slots of its preferred start time $B_i$.
This deviation can either be to the left or to the right of the preferred time slot. We assume here that in terms of
inconvenience, the scheduling of a device $X_i$ time slots before its preferred time slot is equivalent to the inconvenience caused
by scheduling the same device $X_i$ time slots after its preferred time slot.
Since $t \in [1,T]$, therefore if e.g. $B_i=1$ then we can only perform
task $i$ ahead of $B_i$ and schedule it in interval $[B_i,B_i+X_i]$.
Similarly if $B_i=T$ then we can only perform task $i$ before $B_i$ and schedule it in interval $[B_i-X_i, B_i]$.
Thus any non-essential task $i \in \hat{\mathcal{S}}$ willing to tolerate $X$-time slot deviation can be scheduled in the
interval $[\alpha_i,\beta_i]$ where $\alpha_i=\max(1,B_i-X_i)$ and $\beta_i=\min(T-\Delta T_i+1,B_i+X_i)$.
All the non-essential devices $j \notin \hat{\mathcal{S}}$ have $X_j=0$ and they
have to be scheduled exactly at time slot $B_i$ and completed after $\Delta T_j$
time slots. We then treat all such devices $j \notin \hat{\mathcal{S}}$ as essential
load, determine $\bar{S}_j^e=\{S_j(t)\}_{t=1}^{T}$ (as explained in the description of the UC schedule) and then update the essential load accordingly i.e.
\begin{equation}
\tilde{E}(t)= E(t) + \sum_{j \notin \hat{\mathcal{S}}} S_i(t) \:,\forall t
\label{ess_update}
\end{equation}
We can now formulate the scheduling problem as follows (we refer this problem as $\mathcal{P}$),
\begin{equation}
\mathcal{P}: \text{min} \quad \max_{1\leq t \leq T}\quad
\tilde{E}(t)+\sum_{i \in \hat{\mathcal{S}}} S_i(t) \label{sch1_obj}
\end{equation}
subject to,
\begin{equation}
\sum_{t=1}^{T} \sum_{i \in \hat{\mathcal{S}}} S_i(t) = \sum_{i \in \hat{\mathcal{S}}} \tilde{S}_i 
\label{sch1_const1}
\end{equation}
\begin{equation}
\sum_{j=t}^{t+\Delta T_i} S_i(j) = \tilde{S}_i \quad, \forall i \in \hat{\mathcal{S}}, t \in [\alpha_i,\beta_i]
\label{sch1_const2}
\end{equation}
Eq (\ref{sch1_const1}) indicates that the total energy consumed by all the non-essential tasks should be equal to their
total required energy. Eq (\ref{sch1_const2}) says that if non-essential task $i \in \hat{\mathcal{S}}$ starts at time $t$ then it should be finished
at time $t+\Delta T_i$ without interruption. The start time of flexible devices can lie in the interval $t \in [\alpha_i,\beta_i]$.

We can then discuss some special cases of the above general problem.
If all the devices are flexible then $\hat{\mathcal{S}}=\mathcal{S}$
and if all the devices are 100\% flexible then $\alpha_i=1\:,\forall
i$ and $\beta_i=T-\Delta T_i+1 \:,\forall i$ (this value of
$\beta_i$ will allow non-essential task to finish within $T$ time
slots). The solution of this \textbf{100\% flexible problem} is the
best possible practical schedule for the grid and achieves maximum
flatness for given set of essential and non-essential tasks.
Similarly, if all the devices declare themselves as flexible but
allow only $X$-time slot deviation (we assume the same $X$ for all
the devices) then we call this special problem as \textbf{$X$-time
slot deviation problem}. If $Y<K$, devices declare them as 100\%
flexible then we call this special problem as \textbf{$Y$-device
flexible problem}. A \textbf{100\% flexibile device} can be
scheduled at any time $t \in [1, T-\Delta T_i+1]$.

\section{Solution and Algorithm Development}
In this section we discuss the solution of the above scheduling
problems and design practical scheduling algorithms. The optimal
solution of the above problem (including all the special cases) in
general depends on the sequence or order in which we consider
non-essential loads. We illustrate this fact by following simple example.\\
\textbf{Example:} Consider $T=3$ time slots. The essential load is
given as $E(t)=\{2,1,0\}$. There are two 100\% shiftable loads with
demands per time slot given as $\bar{S}_1=\{5,0,0\}$ and
$\bar{S}_2=\{2,2,0\}$. There are two possible permutations, load 1
followed by load 2 or load 2 followed by load 1. In the first case,
the final load per time slot is, \{4,3,5\} with a peak load of 5 in
third time slot. For the second case when load 2 is scheduled before
load 1 we obtain two possible schedules, \{7,3,2\} or \{2,3,7\} both
of which are optimal for this order and give a peak load of 7 in
both schedules. Thus, in order to reduce the peak, we should
schedule load 1 before scheduling load 2. Therefore the sequence in
which we consider non-essential loads cannot be ignored. We now
prove that the above problem $\mathcal{P}$ is NP hard problem.
\begin{Theorem}\label{NP}
The defined problem $\mathcal{P}$ is NP-hard.
\end{Theorem}
\begin{proof}
We consider the special case of the defined problem, where we restrict that $\Delta T_i=1$, and $t\in \{0, T\}, E(t)=0$. We then prove that the special case is NP-hard by a induction from the Multi-Processor Scheduling problem, which is a well-know NP-hard problem in the strong sense.

Multi-Processor Scheduling problem: we are given $m$ identical machines in $\mathcal{M}=\{M_1,M_2,\cdots\}$
and $n$ jobs in $\mathcal{J}=\{J_1,J_2,\cdots, J_n\}$. Job $J_i$ has a processing time $p_i\geq 0$. The objective of Multi-Processor Scheduling problem is to assign jobs to the
machines so as to minimize the maximum load of the machines.

Given an instance of Multi-Processor Scheduling Problem, we can construct an instance of the decision version of the above special case of the defined problem in polynomial time as follows. Let there be $|\mathcal{M}|$ time slots that can be scheduled for tasks, there be $|\mathcal{J}|$ shiftable tasks, and $p_i$ be the power consumed for the $i$-th task. Then, the load of the tasks scheduled at time $t$ is equal to the load of the jobs assigned  at the $t$-th machine. In other words, the objective of minimizing the maximum load at each time is is to minimize the maximum working load assigned at each machine. Thus,
the instance of Multi-Processor Scheduling problem is equivalent to an instance of the special case of the defined load balancing problem. Thus, by induction, the defined load balancing problem is NP-hard.
\end{proof}

Despite the fact that the problems are NP hard we can still design
an algorithm to find the optimal schedules. However the complexity
of the optimal algorithm is exponential which makes it infeasible
when the number of flexible devices is large. We give the optimal
algorithm for problem $\mathcal{P}$ below.

\subsection{Optimal Algorithm}
Let $U=\{(\tilde{S}_1,\Delta T_1,\alpha_1,\beta_1),\ldots,(\tilde{S}_Y,\Delta T_Y) \}$ hold the total power consumption, required completion time and lower and
upper limits of scheduling interval (calculated based on specified $X_i$-time slot deviations) for non-essential tasks willing to be flexible. Note that $Y\leq K$ and $\hat{\mathcal{S}}=\{1,\ldots,Y\}$. Let $\hat{\mathcal{G}}$ denote all possible permutations of this set $\hat{\mathcal{S}}$ i.e. all possible ways of arranging the shiftable tasks in this set. The total number of permutations is $Y!$.
Let $\tilde{L}(t)$ denote the total load scheduled in time slot $t$. \\
\\1. Update essential load according to eq (\ref{ess_update}) for all $j \notin \hat{\mathcal{S}}$.
\\2. For all $m \in \hat{\mathcal{G}}$
\\3. \hspace*{.25 cm}Initialize: $L_m(t)=\tilde{E}(t) \:,\forall t$.
\\4. \hspace*{.25 cm}for all $U_i=(\tilde{S}_i,\Delta T_i, \alpha_i,\beta_i)) \in U$ do,
\\5. \hspace*{.5 cm} for all $t \in [\alpha_i,\beta_i]$
\\6. \hspace*{1 cm} for $j=1,\ldots,\Delta T_i$
\\7. \hspace*{1.5 cm} $X(t,j)=L_m(t+j-1)+\frac{\tilde{S}_i}{\Delta T_i}$
\\8. \hspace*{1 cm} end for
\\9. \hspace*{.5 cm} end for
\\10. \hspace*{.5 cm} $t^*_i=\min_t \max_j X(t,j)$
\\11. \hspace*{.5 cm} $L_m(t^*_i+j-1)=X(t^*_i,j) \:, j=1,\ldots,\Delta T_i$
\\12. \hspace*{.25 cm} end for
\\13. end for
\\14. $m^*=\min_m \max_t L_m(t)$
\\15. $\tilde{L}(t)=L_{m^*}(t)$ \\
\\ In the first line we update the essential load if there are some non-essential devices which declare them as non-flexible.
For each permutation $m$ task $i$ can be scheduled any time between $[\alpha_i,\beta_i]$.
If task $i$ starts at time $t$ then it will be complete at $\Delta T_i-1$. We obtain
all the schedules with all possible start times in lines 5-9 for task $i$. From all
possible schedules we select the one which gives the minimum peak in line 10 and select
this best schedule. In line 11 we update the total load and repeat for the next task $i$.
Finally in line 14 we select the best order $m^*$ in which we should consider the shiftable
tasks. The final schedule is $\tilde{L}(t)$ given in line 15. This is the optimal algorithm.
However, the complexity of this algorithm is exponential which may not be feasible when $Y>>1$.


\subsection{Sub-optimal Algorithm}
We discuss a special case of the above problem when all non-essential tasks have the same power consumption $\frac{\tilde{S}_i}{\Delta T_i}=\hat{S}\:,\forall i \in \hat{\mathcal{S}}$ i.e.
\begin{eqnarray}
S_i(t)=\Bigg\{\begin{array}{c}
                \hat{S},\quad \text{if} \: t=J,\ldots,\Delta T_i-1 \\
                0,\quad \text{otherwise}
                \end{array}
\label{subruly22}
\end{eqnarray}
The required number of time slots to complete each task however are different i.e. $\Delta T_i \neq \Delta T_j$. In this case
the sequence in which we pick the tasks for scheduling becomes irrelevant. Based on this observation we now develop a low complexity sub-optimal algorithm. \\
\\1. Initialize: $\hat{S} = \max_i \frac{\tilde{S}_i}{\Delta T_i} ,\: i \in \hat{\mathcal{S}}$ and $L(t)=\tilde{E}(t) \:,\forall t$.
\\2. \hspace*{.25 cm}for all $U_i=(\hat{S},\Delta T_i,\alpha_i,\beta_i) \in U$ do,
\\3. \hspace*{.5 cm} for all $t \in [\alpha_i,\beta_i]$
\\4. \hspace*{1 cm} for $j=1,\ldots,\Delta T_i$
\\5. \hspace*{1.5 cm} $X(t,j)=L(t+j-1)+\hat{S}$
\\6. \hspace*{1 cm} end for
\\7. \hspace*{.5 cm} end for
\\8. \hspace*{.5 cm} $t^*_i=\min_t \max_j X(t,j)$
\\9. \hspace*{.5 cm} $L(t^*_i+j-1)=X(t^*_i,j) \:, j=1,\ldots,\Delta T_i$
\\10. \hspace*{.25 cm} end for
\\11. for all $i$
\\12. \hspace*{.25 cm} for $t=t^*_i, \ldots, t^*_i+\Delta T_i$
\\13. \hspace*{.5 cm} $L(t)=L(t)-\max\bigg(\hat{S}-\tilde{S}_i (t),0\bigg)$
\\14. \hspace*{.25 cm} end for
\\15. end for

In this algorithm we initially assume that all the shiftable loads have same power consumption $\hat{S}$ per time slot
where $\hat{S}$ is taken as the maximum power consumption across all the non-shiftable devices willing to tolerate $X_i$-time slot deviation.
In lines 2-10, we arbitrarily pick the tasks one after another and find the best scheduling time $t^*_i$ for each shiftable task $i$
in their scheduling interval $t \in [\alpha_i,\beta_i]$. Once we obtain the schedule then in lines 11-15 we restore the loads to their actual power consumption
levels.


\subsection{Comparison of various schedules}
We can measure the difference between any two schedules $\mathcal{R}_n=\{L_n(t)\}_{t=1}^{T}$ and $\mathcal{R}_m=\{L_m(t)\}_{t=1}^{T}$ where $L_i(t)$ denotes the load at time slot $t$ by measuring their mean square error i.e.
\begin{equation}
\mathcal{E}(\mathcal{R}_n,\mathcal{R}_m)=\sum_{t=1}^T \bigg(L_n(t)-L_m(t)\bigg)^2
\label{mse_diff}
\end{equation}
Let $\mathcal{R}n$ denote any arbitrary schedule. As defined before let $\hat{\mathcal{R}}$ denote the GC schedule while
$\tilde{\mathcal{R}}$ denote the UC schedule. Then we define,
\[\gamma=\mathcal{E}(\mathcal{R}_n,\hat{\mathcal{R}})\]
\[\zeta=\mathcal{E}(\mathcal{R}_n,\tilde{\mathcal{R}})\]
where $\gamma$ measures the deviation of any arbitrary schedule
$\mathcal{R}_n$ for the given set of load conditions from the GC
schedule $\hat{\mathcal{R}}$; while $\zeta$ measures the deviation
of any arbitrary schedule $\mathcal{R}_n$ from the UC schedule
$\tilde{\mathcal{R}}$. The smaller the value of $\gamma$ means that
schedule is more flat; while a small value of $\zeta$ means that
schedule is more close to the UC schedule.

\section{Simulation Results}
We consider a generalised simulation setup of residential household
appliances where electricity consumption is assumed to be constant
over the consumption duration and represented in kWh. We generate
essential loads in each time slot as discrete uniform integer random
variables, taking values between 1kWh and 5kWh. Each time slot
represents one hour. In addtion, we assume that there are 100
generalized devices which can be shifted. The total power
consumption of each shiftable device is generated as a discrete
uniformly distributed random variable taking values between 1kWh and
5kWh. The total duration of each shiftable task is generated as a
discrete uniform random integer variable taking values between 1 and
5 time slots. We also assume that each shiftable device has a
preferred time slot. Again this preferred time slot is generated as
a discrete uniform random integer variable.

\begin{figure}[htb]
\centering
\includegraphics[scale=0.4]{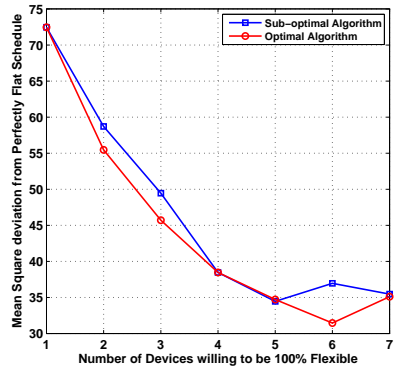}
\vspace{-1em}
\caption{Mean Square Deviation from a Perfectly Flat
Schedule (GC Schedule) ($\gamma$) vs Number of 100\% Shiftable
devices} \label{fig0:fig}
\end{figure}

In Fig. \ref{fig0:fig} we compare the optimal algorithm with the
sub-optimal algorithm. We assume that all the non-essential devices
are 100\% flexible i.e. they can be scheduled at any time $t \in
[1,T-\Delta T_i+1]$. For comparison we measure the mean square
difference of both the schedules from a perfectly flat schedule (GC
Schedule). Since the complexity of the optimal algorithm is
exponential, we restrict ourself to only 7 non-essential devices. We
can see that although there is a small difference between the
performance of both the schedules, the complexity reduction between
the two algorithms is significant, and we will use only the
sub-optimal algorithm in the following simulations.

\vspace{-1em}
\vspace{-.5em}

\begin{figure}[htb]
\begin{minipage} {0.45\linewidth}
\centering
\includegraphics[width=45mm,height=45mm]{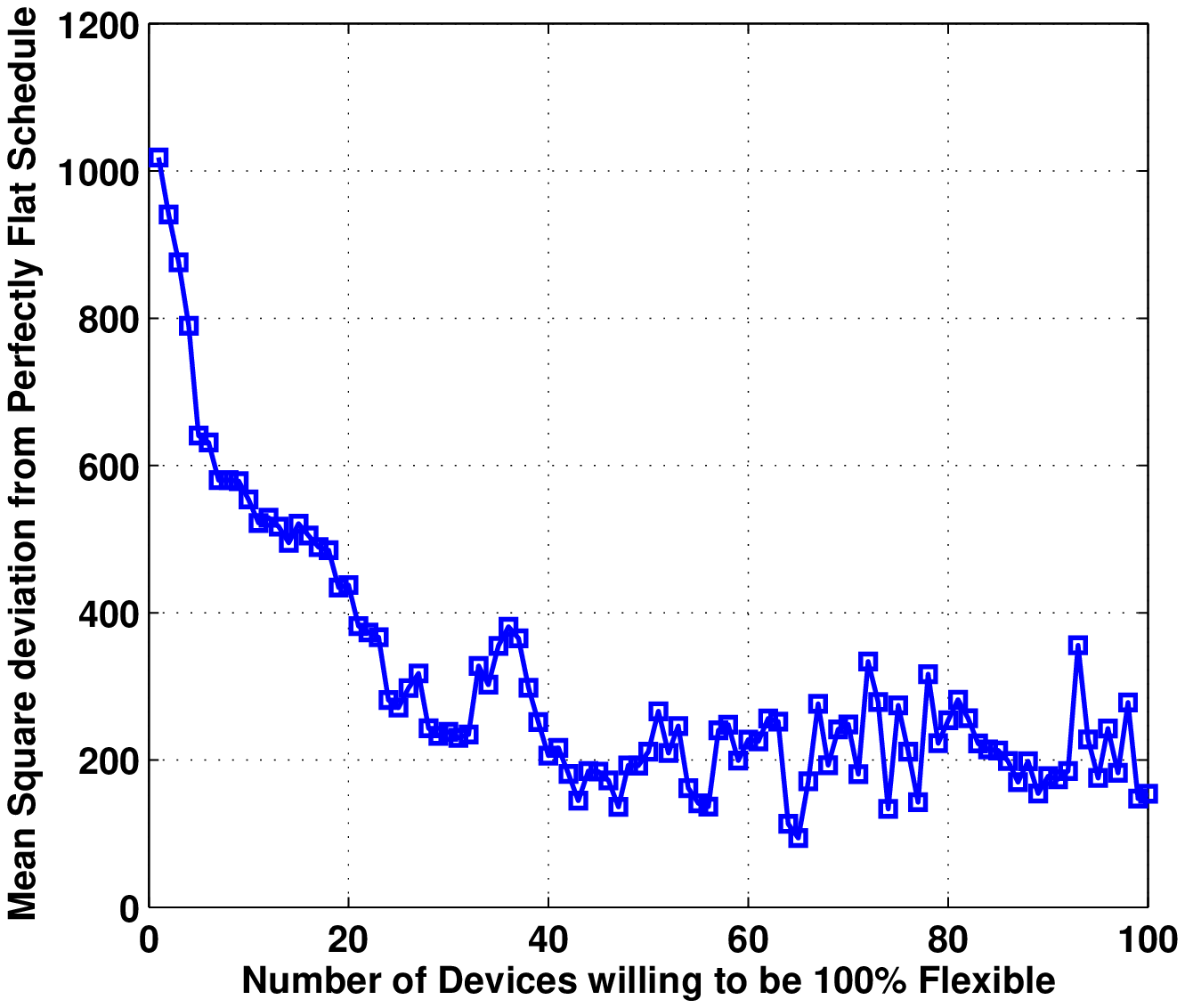}
\caption{Mean Square Deviation from perfectly flat schedule
($\gamma$) vs Number of 100\% Shiftable devices} \label{fig1:fig}
\end{minipage}
\hspace{0.5cm}
\begin{minipage} {0.45\linewidth}
\centering
\includegraphics[width=45mm,height=42mm]{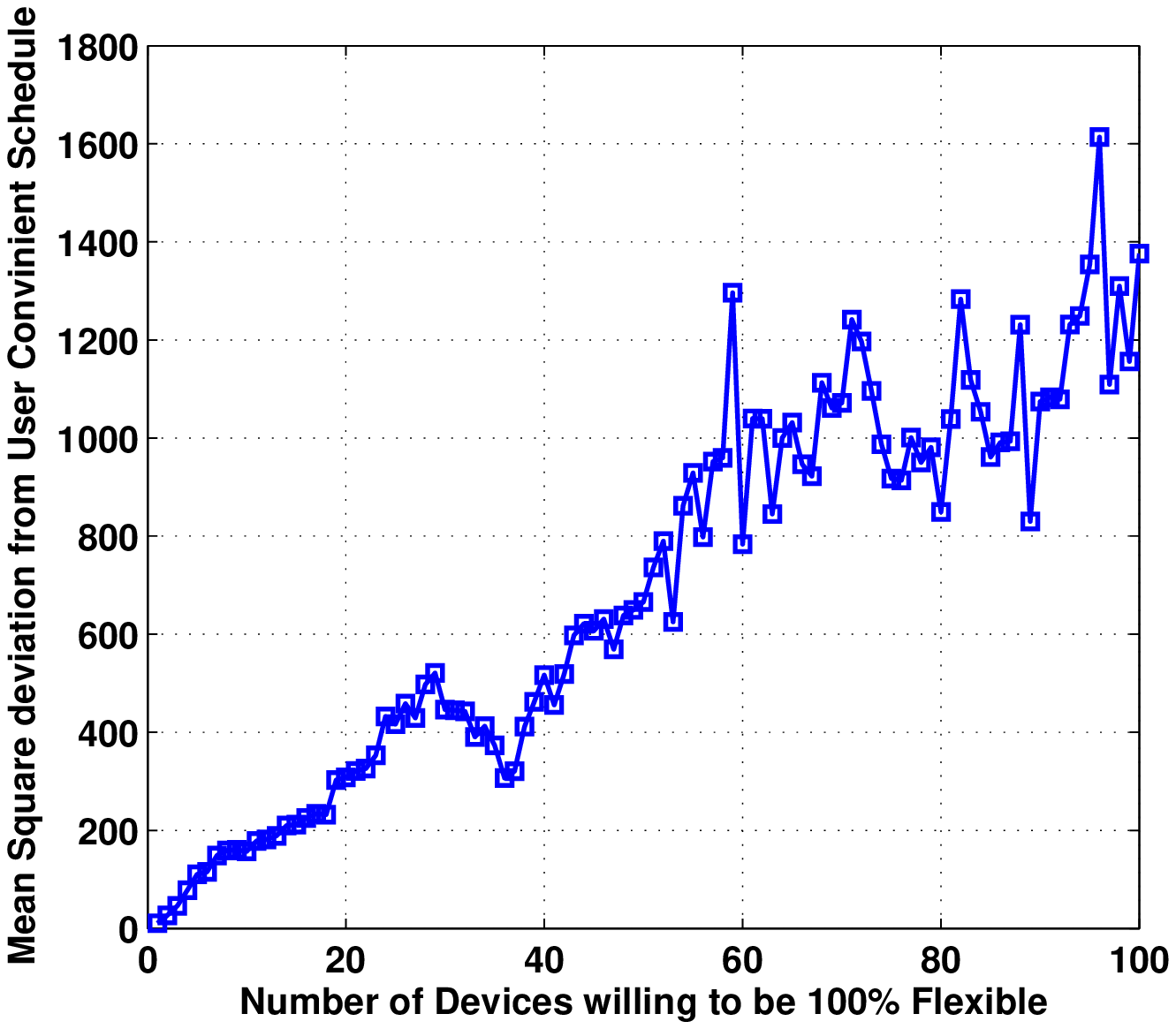}
\caption{Mean Square Deviation from User convenient schedule
($\zeta$) vs Number of 100\% Shiftable devices} \label{fig2:fig}
\end{minipage}
\end{figure}

In Fig. \ref{fig1:fig} and Fig. \ref{fig2:fig} we vary the number of
devices which are willing to be 100\% flexible. All other devices
which are not willing to be flexible are then treated as essential
load and their power consumption is added to the essential load at
their preferred time slots. We plot $\gamma$, the deviation of our
proposed sub-optimal schedule from the GC schedule in Fig.
\ref{fig1:fig}. It is obvious that as more devices become flexible
this deviation decreases. However, we can observe that after 40
devices the value of $\gamma$ does not decrease much which means
that there is not much gain for the grid if more devices become
flexible. The flatness level achieved by 40 devices is comparable to
the flatness level achieved by 100 devices. In Fig. \ref{fig2:fig}
we plot $\zeta$, the deviation of our proposed sub-optimal algorithm
from the UC schedule. As more devices become flexible their
scheduling is not performed at their most convenient time slots and
thus users suffer more inconvenience. The level of inconvenience
keeps on increasing as more devices become flexible. When 40 devices
are 100\% flexible the value of $\zeta$ is 517 and the corresponding
value of $\gamma$ is 206. Similarly when all the devices are 100\%
flexible the value of $\zeta$ is 1375 and that of $\gamma$ is 154.
If we define relative inconvenience level as,
$\hat{\zeta}=\frac{\zeta}{max(\zeta)}\times 100$ and relative
flatness level as $\hat{\gamma}=\frac{\gamma}{max(\gamma)}\times
100$. Then for 40 devices $\hat{\zeta}=32.5$\% and
$\hat{\gamma}=20.2$\% while for 100 devices we have
$\hat{\zeta}=85.2$\% and $\hat{\gamma}=15.1$\%. Thus if a user only
allows 40 devices to become 100\% flexible he can reduce relative
inconviniec level by $85.2-32.5=52.7$\% while the non-flatness will
only increase by $20.2-15.1=5.1$\%. \textbf{This results shows that
there is a minimum level of customer participation in the smart grid
that the grid should aim for that would maximize the gain to the
operator while at the same time imposing minimal inconvenience.
Based on this observation, the inconvenience to the customer will
not be too significant. Although these results may not be
representative of the system, but it does indicate a great research
opportunity to reduce system wide costs at relatively small
individual inconvenience.}

\begin{figure}[htb]
\begin{minipage} {0.45\linewidth}
\centering
\includegraphics[scale=0.32]{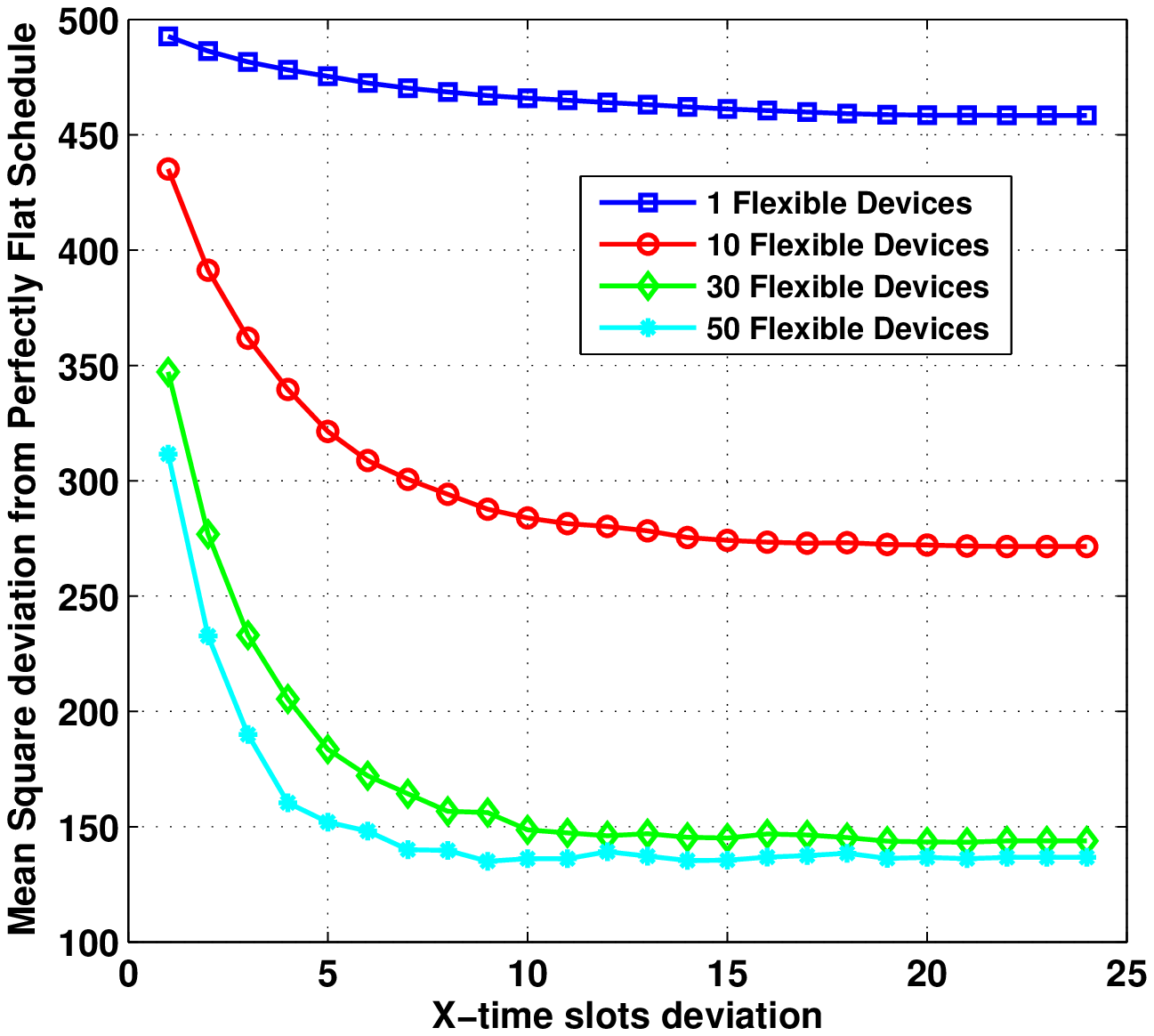}
\caption{Mean Square Deviation from Perfectly Flat Schedule (GC Schedule) ($\gamma$) vs X-time Slot Deviation} \label{fig5:fig}
\end{minipage}
\hspace{0.5cm}
\begin{minipage} {0.45\linewidth}
\centering
\includegraphics[scale=0.32]{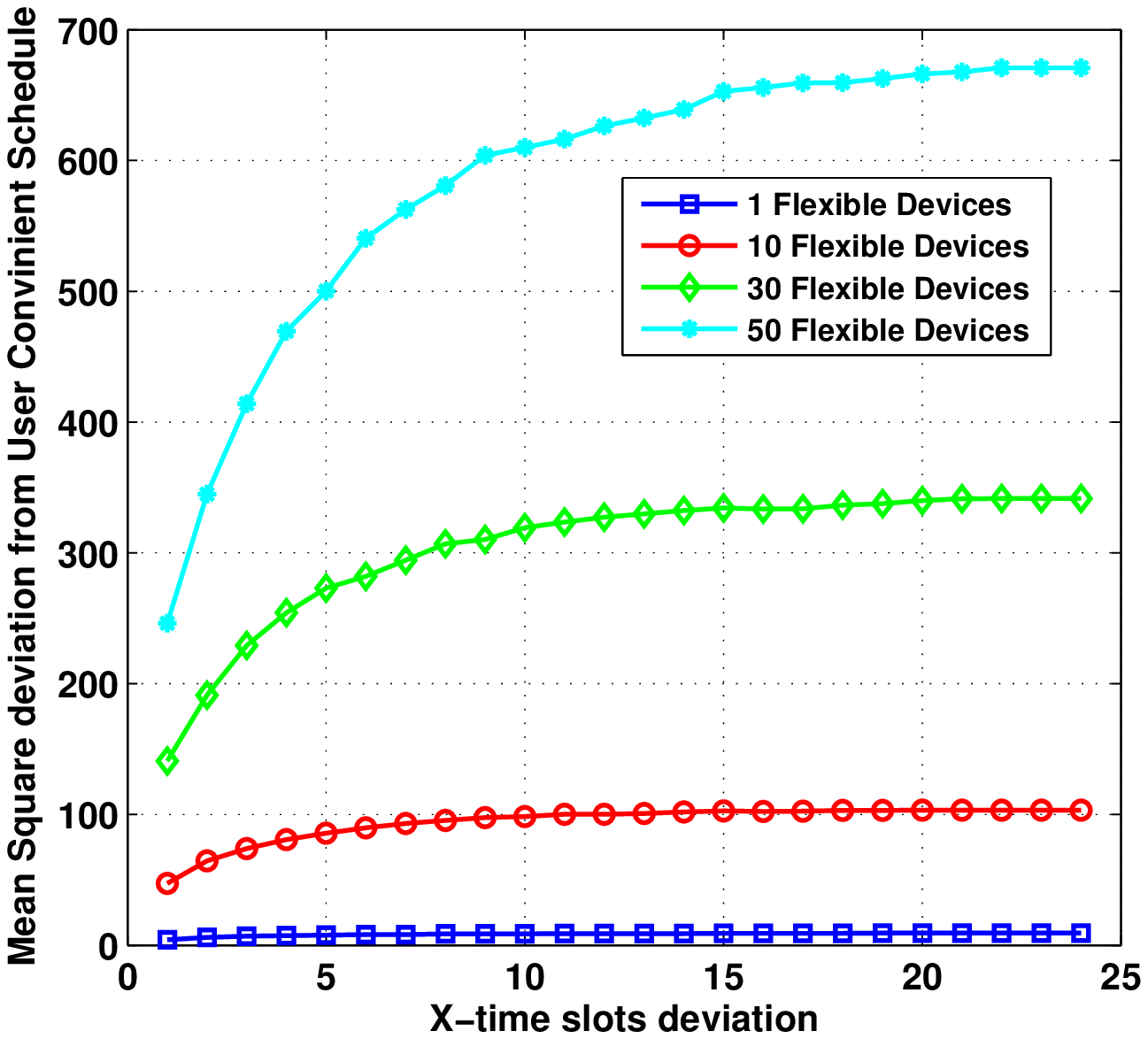}
\caption{Mean Square Deviation from Perfectly Flat Schedule (GC Schedule) ($\gamma$) vs X-time Slot Deviation} \label{fig6:fig}
\end{minipage}
\end{figure}

In Fig. \ref{fig5:fig} and Fig. \ref{fig6:fig}, we obtain various
curves by varying the number of flexible devices. In these
simulations we assume that the total number of available flexible
devices can be up to 50. When the devices declare themselves as
flexible we can schedule them according to their described X-time
slot deviation levels. The load of devices declaring them as
non-flexible is then added to the essential load. E.g. if 10 devices
declare themselves as flexible then the load of remaining 40 devices
is added to the essential load. Therefore the total load in all
these curves is same. As more devices become flexible, we can
achieve more flatness as evident in Fig. \ref{fig5:fig}. However the
gains in flatness diminish and are not very significant as the
number of flexible devices are increased from 30 to 50.\textbf{
Again, this could represent large system savings at minimal
individual costs. }Similarly the gains in flatness also does not
increase much as the X-time slot deviation increases beyond 10 time
slots. On the other hand in Fig. \ref{fig6:fig} we can see that
increasing the number of flexible devices significantly increases
the inconvenience levels of users. When 50 devices are flexible
users experience much more inconvenience compared to 30 flexible
devices. Also increasing X-time slot deviation also increases user
inconvenience. \textbf{Hence, much further research is required to
quantify the trade-off between benefit versus the users
participation and inconvenience caused. }From these observations, in
our test case, we can conclude that significant gains in flatness
can be achieved by declaring a small number of devices as flexible
and keeping X-time slot deviation up to 10 time slots. How this
translates to larger more representative systems need further
examination.

 \vspace{-.5ex}
\section{Conclusion}
In this paper we study the problem of load balancing in smart grids.
We proposed algorithms to obtain schedules by varying the number of
flexible devices and inconvenience levels. Then we identify the
level of compromise between the grid objective of load balancing and
user convenience levels. We show that by allowing only a small
portion of the activities to become flexible, users can contribute
significantly towards load balancing due to aggregating effect.
Similarly by letting the scheduling of activities deviate just a few
hours from their preferred time slots can also significantly impact
load balancing for the grid. More practical system and load models
will be used in the future work to quantify these results. It is
also interesting to investigate what kind of incentives that can be
provided by the grid to encourage the users to have their load be
flexible.

 \vspace{-.5ex}
\section{Acknowledgment}
This research is partly supported by SUTD-ZJU/RES/02/2011,
International Design Center, EIRP, SSE-LUMS via faculty research
startup grant, and Fundamental Research Funds (no. 2012HGBZ0640).

\end{document}